# Some divisibility properties of q-Fibonacci numbers

*Johann Cigler*

Fakultät für Mathematik

Universität Wien

johann.cigler@univie.ac.at

**Abstract**

We give a survey of some known and some new results about factors of different sorts of $q-$Fibonacci numbers.

## 0. Introduction

Let $(F_n)_{n \geq 0} = (0,1,1,2,3,5,8,\cdots)$ be the sequence of Fibonacci numbers and let $v_p(n)$ be the $p-$adic valuation of $n,$ i.e. the highest power of the prime number $p$ which divides $n.$ The Fibonacci numbers satisfy (cf. [7]) $v_5(F_n) = v_5(n),$ $v_2(F_{3n}) = 1$ for odd $n$ and $v_2(F_{6n}) = v_2(n) + 3.$ If $p$ is a prime different from $2$ and $5$ then either $F_{p-1}$ or $F_{p+1}$ is divisible by $p.$

For $q \in \mathbb{C}$ let $[n] = [n]_q = \frac{1-q^n}{1-q} = 1 + q + \cdots + q^{n-1}$ and let $\begin{bmatrix} n \\ k \end{bmatrix} = \frac{[n]\cdots[n-k+1]}{[1][2]\cdots[k]}$ be a $q-$binomial coefficient.

The Schur-Carlitz $q-$Fibonacci numbers $F_n(q) = \sum_{k=0}^{\lfloor \frac{n-1}{2} \rfloor} q^{k^2} \begin{bmatrix} n-1-k \\ k \end{bmatrix}$ and

$G_n(q) = \sum_{k=0}^{\lfloor \frac{n-1}{2} \rfloor} q^{k^2+k} \begin{bmatrix} n-1-k \\ k \end{bmatrix}$ (cf. [9],[2]) which have been introduced by I. Schur in his proof of the Rogers-Ramanujan identities inherit some of the properties for odd primes and the $q-$Fibonacci numbers $f(n,q) = \sum_{k=0}^{\lfloor \frac{n-1}{2} \rfloor} q^{\binom{k}{2}} \begin{bmatrix} n-1-k \\ k \end{bmatrix}$ introduced in [4] inherit divisibility properties by $2.$

## 1. Divisibility properties for odd primes $p \neq 5.$

**1.1.** The $q-$Fibonacci numbers $F_n(q)$ satisfy the recurrence

$$F_n(q) = F_{n-1}(q) + q^{n-2} F_{n-2}(q) \tag{1.1}$$

with initial values $F_0(q) = 0$ and $F_1(q) = 1$.

The first terms are

$$0, 1, 1, 1+q, 1+q+q^2, 1+q+q^2+q^3+q^4, 1+q+q^2+q^3+2q^4+q^5+q^6, \cdots.$$

It is clear that $F_n(1) = F_n$.

**Theorem 1.1** (George E. Andrews, Leonard Carlitz [1])

*If $p$ is an odd prime with $p \equiv \pm 2 \mod 5$ then $F_{p+1}(q) \equiv 0 \mod [p]_q$.*

For $q = 1$ this can be proved (cf. [6], Theorem 180) using Binet's formula

$$F_n = \frac{1}{2^n \sqrt{5}}\left(\left(1+\sqrt{5}\right)^n - \left(1-\sqrt{5}\right)^n\right).$$

This gives

$$2^p F_{p+1} = \binom{p+1}{1} + \binom{p+1}{3} 5 + \cdots + \binom{p+1}{p} 5^{\frac{p-1}{2}}.$$

Here all binomial coefficients are divisible by $p$ except the first and last one. Therefore

$2^p F_{p+1} \equiv 1 + 5^{\frac{p-1}{2}} \mod p$. Hence $F_{p+1} \equiv 0 \mod p$ if $5^{\frac{p-1}{2}} \equiv \left(\frac{5}{p}\right) = -1 \mod p$. By the quadratic reciprocity law $\left(\frac{p}{5}\right)\left(\frac{5}{p}\right) = 1$. This implies $\left(\frac{p}{5}\right) = -1$ and thus $p \equiv \pm 2 \mod 5$.

Since there is no analogue of Binet's formula for $q$–Fibonacci numbers L. Carlitz used the polynomial version of the first Rogers-Ramanujan identity (cf. [9], [5])

$$F_{n+1}(q) = \sum_{k=-\left\lfloor \frac{n+2}{5} \right\rfloor}^{\left\lfloor \frac{n+2}{5} \right\rfloor} (-1)^k q^{\frac{k(5k-1)}{2}} \begin{bmatrix} n \\ \left\lfloor \frac{n+5k}{2} \right\rfloor \end{bmatrix}. \tag{1.2}$$

He showed more generally

**Lemma 1.1**

*Let $\Phi_n(q)$ be the $n$–th cyclotomic polynomial. Then $F_{n+1}(q)$ is divisible by $\Phi_n(q)$ if and only if $n \equiv \pm 2 \mod 5$, where $n$ is an arbitrary positive integer.*

For a prime $n = p$ the cyclotomic polynomial reduces to $1 + q + \cdots + q^{p-1} = [p]_q$ and therefore implies Theorem 1.1.

Since $\begin{bmatrix} n \\ k \end{bmatrix}$ is divisible by $\Phi_n(q)$ for $1 \leq k \leq n-1$ we get by (1.2)

$$F_{n+1}(q) \equiv (-1)^r q^{\frac{r(5r+1)}{2}} \begin{bmatrix} n \\ \frac{n-5r}{2} \end{bmatrix} + (-1)^r q^{\frac{r(5r-1)}{2}} \begin{bmatrix} n \\ \frac{n+5r}{2} \end{bmatrix} \mod \Phi_n(q) \quad (1.3)$$

with $r = \left\lfloor \frac{n+2}{5} \right\rfloor$.

It suffices to verify that $F_{n+1}(q) \equiv 0 \mod \Phi_n(q)$ for $n = 10m+2,\ 10m+3,\ 10m+7,\ 10m+8$.

This is shown in the following table:

| $n$ | $r = \left\lfloor \frac{n+2}{5} \right\rfloor$ | $e(r) = \left\lfloor \frac{n+5r}{2} \right\rfloor$ | $e(-r) = \left\lfloor \frac{n-5r}{2} \right\rfloor$ | $\begin{bmatrix} n \\ e(r) \end{bmatrix}$ | $\begin{bmatrix} n \\ e(-r) \end{bmatrix}$ |
|---|---|---|---|---|---|
| $10m+2$ | $2m$ | $10m+1$ | $1$ | $[n]$ | $[n]$ |
| $10m+3$ | $2m+1$ | $10m+4$ | $-1$ | $0$ | $0$ |
| $10m+7$ | $2m+1$ | $10m+6$ | $1$ | $[n]$ | $[n]$ |
| $10m+8$ | $2m+2$ | $10m+9$ | $-1$ | $0$ | $0$ |

In each case both terms of (1.3) vanish modulo $\Phi_n(q)$.

Also observe that $f(q) \equiv 0 \mod \Phi_n(q)$ for a polynomial $f(q)$ is equivalent with $f(\zeta_n) = 0$ for a primitive $n$-th root of unity $\zeta_n$.

**1.2.** The $q$-Fibonacci numbers $G_n(q)$ satisfy the recurrence

$$G_n(q) = G_{n-1}(q) + q^{n-1} G_{n-2}(q) \quad (1.4)$$

with initial values $G_0(q) = 0$ and $G_1(q) = 1$. The first terms are

$0, 1, 1, 1+q^2, 1+q^2+q^3, 1+q^2+q^3+q^4+q^6, 1+q^2+q^3+q^4+q^5+q^6+q^7+q^8, \cdots$.

The polynomial version of the second Rogers-Ramanujan identity (cf.[9],[5]) gives

$$G_n(q) = \sum_{k=-\left\lfloor \frac{n+2}{5} \right\rfloor}^{\left\lfloor \frac{n+2}{5} \right\rfloor} (-1)^k q^{\frac{k(5k-3)}{2}} \begin{bmatrix} n \\ \frac{n-1+5k}{2} \end{bmatrix}. \quad (1.5)$$

For $n = 5m$ this implies

$$G_{5m}(q) = \sum_{k=-m}^{m} (-1)^k q^{\frac{k(5k-3)}{2}} \left[ \begin{array}{c} 5m \\ \left\lfloor \frac{5(m+k)-1}{2} \right\rfloor \end{array} \right]_q \equiv 0 \bmod \Phi_{5m}(q) \tag{1.6}$$

since no $q-$binomial coefficient reduces to 1.

As has been observed by H. Pan [8] for $n \not\equiv 0 \bmod 5$ there remain modulo $\Phi_n(q)$ only the terms with $k = r(n)$, where $r(n) = \left\lfloor \frac{n+2}{5} \right\rfloor$ if $n \equiv 3 \bmod 5$ and $n \equiv 4 \bmod 5$ and

$r(n) = -\left\lfloor \frac{n+2}{5} \right\rfloor$ if $n \equiv 1 \bmod 5$ or $n \equiv 2 \bmod 5$.

This leads to the following table where the congruences are modulo $\Phi_n(q)$.

| $n$ | $r(n)$ | $G_n(q)$ | |
|---|---|---|---|
| $5m$ | $0$ | $0$ | |
| $5m+1$ | $-m$ | $(-1)^m q^{\frac{m(5m+3)}{2}}$ | $\equiv q^m$ |
| $5m+2$ | $-m$ | $(-1)^m q^{\frac{m(5m+3)}{2}}$ | $\equiv -q^{3m+1}$ |
| $5m+3$ | $m+1$ | $(-1)^{m+1} q^{\frac{(m+1)(5m+2)}{2}}$ | $\equiv -q^{2m+1}$ |
| $5m+4$ | $m+1$ | $(-1)^{m+1} q^{\frac{(m+1)(5m+2)}{2}}$ | $\equiv q^{4m+3}$ |

The congruences in the right column are easily verified. For example we have for $n = 5m+2$ and even $m$

$$(-1)^m q^{\frac{m(5m+3)}{2}} \equiv q^{\frac{m}{2}(5m+2)} q^{\frac{m}{2}} \equiv -q^{\frac{m}{2}+\frac{5m+2}{2}} = -q^{3m+1}$$

and for odd $m$

$$(-1)^{m+1} q^{\frac{m(5m+3)}{2}} \equiv q^{\frac{m(5m+3)}{2} - \frac{(5m+2)(m-1)}{2}} \equiv q^{3m+1}.$$

**Theorem 1.2** ( H. Pan [8])

*If $p$ is a prime with $p \equiv \pm 1 \bmod 5$ then $G_{p-1}(q) \equiv 0 \bmod [p]_q$.*

For example

$G_{10}(q) = [11]_q [5]_{q^2} \left(1 - q + q^3 - q^4 + q^6\right).$

Let me sketch H. Pan's proof.

By (1.5) we get

$$G_{n-1}(q) = \sum_{k=-\lfloor\frac{n+1}{5}\rfloor}^{\lfloor\frac{n+1}{5}\rfloor} (-1)^k q^{\frac{k(5k-3)}{2}} \begin{bmatrix} n-1 \\ \lfloor\frac{n-2+5k}{2}\rfloor \end{bmatrix}.$$

For $n = 5m+1$ this reduces to

$$G_{5m}(q) = \sum_{k=-m+1}^{m} (-1)^k q^{\frac{k(5k-3)}{2}} \begin{bmatrix} 5m \\ \lfloor\frac{5(m+k)-1}{2}\rfloor \end{bmatrix} = \sum_{k=-m+1}^{m} (-1)^k q^{\frac{k(5k-3)}{2}} \prod_{j=1}^{\lfloor\frac{5(m+k)-1}{2}\rfloor} \frac{[5m+1-j]_q}{[j]_q}$$

$$= \sum_{k=-m+1}^{m} (-1)^k q^{\frac{k(5k-3)}{2}} \prod_{j=1}^{\lfloor\frac{5(m+k)-1}{2}\rfloor} \frac{[5m+1]_q - [j]_q}{q^j [j]_q} \equiv \sum_{k=-m+1}^{m} (-1)^{k+\lfloor\frac{5(m+k)-1}{2}\rfloor} q^{\frac{k(5k-3)}{2} - \binom{\lfloor\frac{5(m+k)+1}{2}\rfloor}{2}} \mod \Phi_n(q).$$

Now observe that

$$\ell(m,k) = \frac{k(5k-3)}{2} - \binom{\lfloor\frac{5(m+k)+1}{2}\rfloor}{2} \quad \text{satisfies} \quad \ell(m, 2k-1) - \ell(m, 2k) = 5m+1 = n$$

if $m \equiv 0 \mod 2$ and $\ell(m, 2k+1) - \ell(m, 2k) = -5m-1 = -n$ if $m \equiv 1 \mod 2$.

Therefore each pair of adjacent terms in $G_m(q) = \sum_{k=-m+1}^{m} (-1)^{k+\lfloor\frac{5(m+k)-1}{2}\rfloor} q^{\ell(m,k)} \mod \Phi_n(q)$

satisfies $\pm q^{\ell(m,2k-1)} \mp q^{\ell(m,2k)} = 0 \mod \Phi_n(q)$ if $m$ is even and
$\pm q^{\ell(m,2k)} \mp q^{\ell(m,2k+1)} = 0 \mod \Phi_n(q)$ if $m$ is odd.

For $n = 5m+4$ and

$$G_{5m+3}(q) = \sum_{k=-m}^{m+1} (-1)^k q^{\frac{k(5k-3)}{2}} \begin{bmatrix} 5m+3 \\ \lfloor\frac{5(m+k)+2}{2}\rfloor \end{bmatrix} \mod \Phi_n(q)$$

the situation is analogous.

With the same arguments H. Pan has shown that

$$F_{5n}(q) \equiv 0 \mod \Phi_{5n}(q). \tag{1.7}$$

By (1.2) we get

$$F_{5n}(q) = \sum_{k=-n+1}^{n} (-1)^k q^{\frac{k(5k-1)}{2}} \begin{bmatrix} 5n-1 \\ \frac{5n+5k-1}{2} \end{bmatrix}_q$$

and as above each pair of adjacent elements sums to 0.

**1.3.** Let $A(x) = \begin{pmatrix} 1 & x \\ 1 & 0 \end{pmatrix}$. Then it is easily verified (cf. [3]) that

$$A(q^{n-1})A(q^{n-2})\cdots A(q)A(1) = \begin{pmatrix} F_{n+1}(q) & G_n(q) \\ F_n(q) & G_{n-1}(q) \end{pmatrix}. \tag{1.8}$$

If we take the determinant of (1.8) we get the $q-$ Cassini formula

$$F_{n+1}(q)G_{n-1}(q) - F_n(q)G_n(q) = (-1)^n q^{\binom{n}{2}}. \tag{1.9}$$

If $q$ is a primitive $n-$ th root of unity then $q^{\binom{n}{2}} = \left(q^{\frac{n}{2}}\right)^{n-1} = -1$ if $n \equiv 0 \bmod 2$ and

$q^{\binom{n}{2}} = (q^n)^{\frac{n-1}{2}} = 1$ if $n \equiv 1 \bmod 2$. Therefore we get

$$(-1)^n q^{\binom{n}{2}} \equiv -1 \bmod \Phi_n(q). \tag{1.10}$$

The above results and Cassini's formula give

**Corollary 1.1**

If $n \equiv 0 \bmod 5$ then $F_n(q)G_n(q) \equiv 0 \bmod \Phi_n(q)$ and if

$n \not\equiv 0 \bmod 5$ then

$$F_n(q)G_n(q) \equiv 1 \bmod \Phi_n(q). \tag{1.11}$$

More generally we get

**Corollary 1.2**

Let $\zeta_k$ be a primitive $k-$ th root of unity. Then

$$\begin{aligned} F_{kn}(\zeta_k) &= F_n F_k(\zeta_k), \\ G_{kn}(\zeta_k) &= F_n G_k(\zeta_k), \end{aligned} \tag{1.12}$$

and therefore

$$F_{kn}(\zeta_k)G_{kn}(\zeta_k) = \begin{cases} 0 & \text{if } k \equiv 0 \bmod 5 \\ F_n^2 & \text{if } k \not\equiv 0 \bmod 5. \end{cases} \tag{1.13}$$

**Proof**

Let $j = mk + \ell$ with $0 \le \ell < k$. Then for $km + \ell \le k(n-m) - \ell - 1$

$$\begin{bmatrix} k(n-m)-\ell-1 \\ km+\ell \end{bmatrix}_q = \prod_{i=1}^{m} \frac{1-q^{k(n-m-i)}}{1-q^{ki}} * \frac{\left(1-q^{k(n-m)-\ell-1}\right)\cdots\left(1-q^{k(n-m)-k+1}\right)\left(1-q^{k(n-m)-k-1}\right)\cdots\left(1-q^{k(n-m)-k-\ell}\right)}{(1-q)\cdots(1-q^{k-1})}$$

$$*\cdots* \frac{\left(1-q^{k-\ell-1}\right)\cdots\left(1-q^{k-2\ell-1}\right)}{(1-q)\cdots(1-q^{\ell})}$$

If we let $q \to \zeta_k$ then the first term converges to $\binom{n-m-1}{m}$, the middle terms give 1 because the factors of the numerator are a permutation of the factors of the denominator, and the last term converges to $\begin{bmatrix} k-\ell-1 \\ \ell \end{bmatrix}_{\zeta_k}$.

Therefore we get

$$F_{kn}(\zeta_k) = \sum_j \begin{bmatrix} kn-1-j \\ j \end{bmatrix}_{\zeta_k} \zeta_k^{j^2} = \sum_m \sum_\ell \begin{bmatrix} k(n-m)-1-\ell \\ km+\ell \end{bmatrix}_{\zeta_k} \zeta_k^{(km+\ell)^2}$$

$$= \sum_m \sum_\ell \binom{n-m-1}{m} \begin{bmatrix} k-\ell-1 \\ \ell \end{bmatrix}_{\zeta_k} \zeta_k^{\ell^2} = \sum_\ell \begin{bmatrix} k-\ell-1 \\ \ell \end{bmatrix}_{\zeta_k} \zeta_k^{\ell^2} \sum_m \binom{n-m-1}{m} = F_n F_k(\zeta_k).$$

The proof for $G_n(q)$ is essentially the same.

**2. The main result for $F_{5n}(q)$ and $G_{5n}(q)$.**

**Theorem 2.1**

Let $n = 5^k m$ with $k \ge 1$ and $m \not\equiv 0 \bmod 5$. Then

$$F_{5^k m}(q) \text{ and } G_{5^k m}(q) \text{ are divisible by } [5^k]_{q^m}. \qquad (2.1)$$

For example

$F_5(q) = [5]_q$,

$F_{10}(q) = [5]_{q^2}\left(1+q+q^4[9]_q\right) = [5]_q\left(1-q+q^2-q^3+q^4\right)\left(1+q+q^4[9]_q\right)$,

$G_5(q) = [5]_q(1-q+q^2)$,

$G_{10}(q) = [5]_{q^2}[11]_q\left(1-q+q^3-q^4+q^6\right)$.

Let us first recall how to prove that $v_5(F_n) = v_5(n)$. By Binet's formula we get

$$F_n = \frac{1}{2^n\sqrt{5}}\left((1+\sqrt{5})^n - (1-\sqrt{5})^n\right) = \frac{1}{2^{n-1}}\sum_{k=0}^{n}\binom{n}{2k+1}5^k = \frac{1}{2^{n-1}}\sum_{k=0}^{n}\frac{n}{2k+1}\binom{n-1}{2k}5^k.$$

For each $k > 0$ we have $v_5\left(\frac{5^k n}{2k+1}\right) > v_5(n)$ and for $k = 0$ we have $v_5\left(\frac{n}{1}\binom{n-1}{0}\right) = v_5(n)$.

This implies $v_5(F_n) = v_5(n)$.

It is rather trivial that $F_{5n}(q)$ and $G_{5n}(q)$ are divisible by $[5]_q$.

To show this observe that $q^n \equiv q^{n(\mathrm{mod}\,5)} (\mathrm{mod}\,[5]_q)$. Therefore $F_{5n}(q) \equiv 0\,(\mathrm{mod}\,[5]_q)$ by (1.1) implies $F_{5n+2}(q) \equiv F_{5n+1}(q)$, $F_{5n+3}(q) \equiv F_{5n+2}(q) + qF_{5n+1}(q) \equiv F_{5n+1}(q)(1+q)$,

$F_{5n+4}(q) \equiv F_{5n+3}(q) + q^2 F_{5n+2}(q) \equiv (1+q+q^2)F_{5n+1}(q)$ and finally

$F_{5n+5}(q) \equiv F_{5n+4}(q) + q^3 F_{5n+3}(q) \equiv (1+q+q^2+q^3+q^4)F_{5n+1}(q) \equiv 0\,(\mathrm{mod}\,[5]_q)$.

Analogously $G_{5n}(q) \equiv 0\,(\mathrm{mod}\,[5]_q)$ by (1.1) implies $G_{5n+2}(q) \equiv G_{5n+1}(q)$,

$G_{5n+3}(q) \equiv G_{5n+2}(q) + q^2 G_{5n+1}(q) \equiv G_{5n+1}(q)(1+q^2)$,

$G_{5n+4}(q) \equiv G_{5n+3}(q) + q^3 G_{5n+2}(q) \equiv (1+q^2+q^3)G_{5n+1}(q)$ and finally

$G_{5n+5}(q) \equiv G_{5n+4}(q) + q^4 G_{5n+3}(q) \equiv (1+q^2+q^3+q^4+q^6)G_{5n+1}(q) = [5]_q(1-q+q^2)G_{5n+1}(q)$

$\equiv 0\,(\mathrm{mod}\,[5]_q)$.

For the general case observe that by (1.7) and (1.6) $F_{5^\ell r}(q) \equiv 0 \mod \Phi_{5^\ell r}(q)$ and $G_{5^\ell r}(q) \equiv 0 \mod \Phi_{5^\ell r}(q)$ for each factor $5^\ell r$ of $5^k m$ with $\ell \geq 1$ and that all $\Phi_{5^\ell r}(q)$ are irreducible. Therefore the product of all these cyclotomic polynomials divides $F_{5^k m}(q)$ and $G_{5^k m}(q)$. But this product coincides with $[5^k]_{q^m}$ because $[5^k]_{q^m} = \frac{1-q^{5^k m}}{1-q^m} = \frac{\prod_{d|5^k m}\Phi_d(q)}{\prod_{d|m}\Phi_d(q)}$.

For example we see that $F_{10}(q)$ is divisible by $\Phi_5(q)\Phi_{10}(q)$, $F_{15}(q)$ is divisible by $\Phi_5(q)\Phi_{15}(q)$, or $F_{20}(q)$ is divisible by $\Phi_5(q)\Phi_{10}(q)\Phi_{20}(q)$.

### 3. The Fibonacci numbers $f_r(n,q)$.

Let for some $r \in \mathbb{Z}$

$$f_r(n,q) = \sum_{k=0}^{\lfloor\frac{n-1}{2}\rfloor} q^{\binom{k}{2}+2rk}\begin{bmatrix}n-1-k\\k\end{bmatrix}. \tag{3.1}$$

These polynomials satisfy the recurrence (cf. [4])

$$f_r(n,q) = f_r(n-1,q) + q^{n-3+2r} f_r(n-3,q) + q^{n-4+4r} f_r(n-4,q) \tag{3.2}$$

with initial values $f_r(0,q) = 0$, $f_r(1,q) = 1$, $f_r(2,q) = 1$, $f_r(3,q) = 1+q^{2r}$ and $f_r(4,q) = 1+q^{2r}+q^{1+2r}$.

Of special interest is $f(n,q) = f_0(n,q)$. The first terms of $f(n,q)$ are

$$0, 1, 1, 2, 2+q, 2+2q+q^2, 2(1+q)(1+q^2), 2+2q+2q^2+4q^3+2q^4+q^5, \cdots.$$

**Conjecture 3.1**

Let $n = 2^k(2m+1)$ with $k \geq 0$. Then

$$f(6n,q) = f\left(6 \cdot (2m+1) \cdot 2^k, q\right) \text{ is divisible by } 2\left[2^{k+2}\right]_{q^{2m+1}}. \tag{3.3}$$

For example $f(12,q) = 2[8]_q \left(1+q^3+q^5+q^6+q^7+2q^8+q^9+q^{11}\right)$ and $f(18,q)$ is divisible by $2[4]_{q^3}$.

Let me prove some trivial facts:

The Fibonacci numbers $F_n$ satisfy $F_{6n} \equiv 0 \mod 8$. For
$(F_n \mod 8)_{n \geq 0} = (0,1,1,2,3,5,0,1,1,2,3,5,0,\cdots).$

Now observe that $f(3n,q)$ is even because $\sum_{k=0}^{\lfloor \frac{3n-1}{2} \rfloor} (-1)^k q^{\binom{k}{2}} \begin{bmatrix} 3n-1-k \\ k \end{bmatrix} = 0.$

Cf. [4], Theorem 3.2 and the literature cited there.

It is also easy to show by induction that

$f(3n,q) \mod 2 = 0$, $f(3n+1,q) \mod 2 = q^{\frac{n(3n-1)}{2}}$, $f(3n+2,q) \mod 2 = q^{\frac{n(3n+1)}{2}}$.

Observe that $f(6,q) = 2(1+q+q^2+q^3) = 2[4]_q$.

The sequence $(q^n \mod [4]) = (1,q,q^2,-1-q-q^2,\cdots)$ is periodic with period 4.

This implies that the sequence $f(n+24,q) \mod [4]_q$ satisfies the same recurrence. It is easily verified that it also has the same initial values. Therefore the sequence $f(n,q) \mod [4]_q$ has period 24. Since it satisfies $f(6n,q) \equiv 0 \mod [4]_q$ we finally get that $f(6n,q)$ is divisible by $2(1+q)(1+q^2)$.

For general $r$ we get

**Conjecture 3.2**

Let $n = 2^k(2m+1)$ with $k \geq 0$. Then

$f_r(6n,q) = f\left(6 \cdot (2m+1) \cdot 2^k, q\right)$ is divisible by $\left[2^{k+2}\right]_{q^{2m+1}}$.

For example $f_r(6,q) = (1+q^{2r})(1+q^{2r+1}+q^{2r+2}+q^{2r+3})$ is a multiple of $(1+q)(1+q^2)$. If $r$ is even then $i$ and $-1$ are roots of the second factor, if $r$ is odd then $i$ is a root of the first factor and $-1$ is a root of the second factor of $f_r(6,q)$.